\documentclass[conference]{IEEEtran}
\IEEEoverridecommandlockouts
\usepackage{cite}
\usepackage{amsmath,amssymb,amsfonts}
\usepackage{graphicx}
\usepackage{textcomp}
\usepackage{xcolor}
\usepackage[caption=false,font=normalsize,labelfon
t=sf,textfont=sf]{subfig}
\usepackage{algorithmic}
\usepackage[ruled,norelsize,vlined,linesnumbered]{algorithm2e}

\def\BibTeX{{\rm B\kern-.05em{\sc i\kern-.025em b}\kern-.08em
    T\kern-.1667em\lower.7ex\hbox{E}\kern-.125emX}}
\begin{document}

\title{
On the Effects of Smoothing Rugged Landscape by Different Toy Problems: A Case Study on UBQP}

\author{\IEEEauthorblockN{Anonymous Authors}}
\author{\IEEEauthorblockN{1\textsuperscript{st} Wei Wang}
\IEEEauthorblockA{\textit{School of Mathematics and Statistics} \\
\textit{Xi'an Jiaotong University}\\
Xi'an, China \\
ww123@stu.xjtu.edu.cn}
\and
\IEEEauthorblockN{2\textsuperscript{nd} Jialong Shi}
\IEEEauthorblockA{\textit{School of Mathematics and Statistics} \\
\textit{Xi'an Jiaotong University}, Xi'an, China \\
\textit{Sichuan Digital Economy Industry }\\ \textit{Development Research Institute}, Chengdu, China \\
jialong.shi@xjtu.edu.cn}
\and
\IEEEauthorblockN{3\textsuperscript{rd} Jianyong Sun}
\IEEEauthorblockA{\textit{School of Mathematics and Statistics} \\
\textit{Xi'an Jiaotong University}\\
Xi'an, China \\
jy.sun@xjtu.edu.cn}
\and
\IEEEauthorblockN{4\textsuperscript{th} Arnaud Liefooghe}
\IEEEauthorblockA{\textit{LISIC}\\
\textit{Université du Littoral}\\\textit{Côte d'Opale (ULCO)}\\
Calais, France\\
arnaud.liefooghe@univ-littoral.fr}
\and
\IEEEauthorblockN{5\textsuperscript{th} Qingfu Zhang}
\IEEEauthorblockA{\textit{Department of Computer Science} \\
\textit{City University of Hong Kong}\\
Hong Kong, China \\
qingfu.zhang@cityu.edu.hk}
\and
\IEEEauthorblockN{6\textsuperscript{th} Ye Fan}
\IEEEauthorblockA{\textit{School of Electronics and Information} \\
\textit{Northwestern Polytechnical University},
Xi'an, China \\
\textit{Research and Development Institute of}\\ \textit{Northwestern Polytechnical University in Shenzhen}\\ Shenzhen, China\\
fanye@nwpu.edu.cn}
}

\maketitle

\begin{abstract}
The hardness of the Unconstrained Binary Quadratic Program (UBQP) problem is due its rugged landscape. Various algorithms have been proposed for UBQP, including the Landscape Smoothing Iterated Local Search (LSILS). Different from other UBQP algorithms, LSILS tries to smooth the rugged landscape by building a convex combination of the original UBQP and a toy UBQP. In this paper, our study further investigates the impact of smoothing rugged landscapes using different toy UBQP problems, including a toy UBQP with matrix $\boldsymbol{\hat{Q}}^1$ (construct by ``+/-1''), a toy UBQP with matrix $\boldsymbol{\hat{Q}}^2$ (construct by ``+/-i'') and a toy UBQP with matrix $\boldsymbol{\hat{Q}}^3$ (construct randomly). We first assess the landscape flatness of the three toy UBQPs. Subsequently, we test the efficiency of LSILS  with different toy UBQPs. Results reveal that the toy UBQP with $\boldsymbol{\hat{Q}}^1$ (construct by ``+/-1'') exhibits the flattest landscape among the three, while the toy UBQP with $\boldsymbol{\hat{Q}}^3$ (construct randomly) presents the most non-flat landscape. Notably, LSILS using the toy UBQP with $\boldsymbol{\hat{Q}}^2$ (construct by ``+/-i'') emerges as the most effective, while $\boldsymbol{\hat{Q}}^3$ (construct randomly) has the poorest result. These findings contribute to a detailed understanding of landscape smoothing techniques in optimizing UBQP.
\end{abstract}


\section{Introduction}
The Unconstrained Binary Quadratic Program (UBQP) \cite{UBQP} has attracted much attention in the field of combinatorial optimization due to its practical applications in financial analysis~\cite{financial}, molecular conformation~\cite{molecular}, traffic management~\cite{traffic} and so on. The aim of this NP-hard problem is to determine binomial decision variables that will maximize a quadratic objective function by choosing binary decision variables. The UBQP problem can be formalized as follows \cite{UBQPformulation}:
\begin{equation}
	\begin{aligned}
		&\textrm{maximize} \quad f(\boldsymbol{x}) = \boldsymbol{x}^T\boldsymbol{Q}\boldsymbol{x} = \sum_{i=1}^{n}\sum_{j=1}^{n}Q_{ij}x_ix_j \\
		&\textrm{subject to}  \qquad \qquad \boldsymbol{x}\in \{0,1\}^n
	\end{aligned}
	\label{ubqp}
\end{equation}
where $\boldsymbol{Q} = [Q_{ij}]$ is an $n \times n$-dimensional matrix and $\boldsymbol{x}$ is an $n$-dimensional vector with binary variables. 
UBQP serves as a common model for a wide range of combinatorial optimization problems, such as maximum cut problems \cite{maxcut} and set partitioning problems \cite{setpartition}, among many others~\cite{otherproblem}.

There exist difficulties in finding a global optimal solution within polynomial time because of the NP-hardness of UBQP~\cite{np_hard}. Numerous heuristics and metaheuristics have been developed over the years in order to cope with this computational challenge. UBQP is characterized by a rugged and irregular fitness landscape, which makes the task of solving it challenging, particularly when there are numerous local optima within the search landscape \cite{UBQPworst}. However, to the best of our knowledge, very few research efforts have been committed to smoothing the landscape of UBQP.

In an effort, Wang et al. \cite{HCUBQP} introduce a landscape smoothing method called Homotopic Convex (HC) transformation for the UBQP by employing a toy UBQP problem formed by $1$s and $-1$s. They design the Landscape Smoothing Iterated Local Search~(LSILS) algorithm based on this HC transformation. However, no study has been conducted on the smoothing effects of HC transformation using different toy problems. This paper extends the investigation into the smoothing effects of HC transformation by introducing two additional toy UBQPs, where the assumption is that toy problems with different landscape flatness have different landscape smoothing effects. Subsequently, employing the previously proposed LSILS algorithm, we evaluate the performance of it using the three toy UBQPs and compare it against a baseline ILS \cite{ILS}. Results indeed indicate that different toy UBOPs do have different effects on smoothing the landscape.

The rest of this paper is organized as follows. Section \ref{2} reviews the related work. In Section \ref{3}, we briefly review the LSILS presented in \cite{HCUBQP} and introduce two LSILS algorithms with different toy problems. Detailed experimental results are reported in Section \ref{4}. Section \ref{5} presents the conclusion.

\section{Related Work}\label{2}
Generally speaking, there exist two types of algorithms for UBQP: exact methods and heuristics. Since the UBQP problem is NP-hard~\cite{UBQP,np_hard}, exact algorithms tend to be time-consuming for large-scale problems. They are not able to obtain the global optimum within limited computing time. Different from exact algorithms, heuristics can obtain near-optimal solutions within a reasonable time for large-scale UBQP instances. Among those are algorithms based on tabu search \cite{DDTS}, \cite{TS}, simulated annealing \cite{SA1}\cite{SA2}, iterated local search \cite{ILS1}, evolutionary algorithm (EA) \cite{EA}, among others. However, we find that the previous algorithms are rarely based on landscape smoothing. In \cite{HCUBQP}, Wang et al. propose a landscape smoothing method HC transformation for UBQPs and then design a LSILS algorithm based on this method.

The behavior of heuristics is significantly influenced by the characteristics of fitness landscapes. Therefore, nowadays, the characteristics of fitness landscapes have been widely studied in many combinatorial optimization problems~\cite{fitness_landscape}. The difficulty to explore a landscape lies in its ruggedness which is directly related to the number of local optima in the landscape. A great amount of efforts have been made on fitness landscape analysis. With regard to UBQP, Tari et al. found a big-valley structure in UBQP landscapes \cite{SampledWalk}\cite{UBQPworst}. They point out that local optima tend to be grouped in a sub-region of the landscape. Merz and Katayama \cite{MA} conducted a landscape analysis of UBQP and observed that local optima are contained in a small fraction of the search space, which corroborates previous research. Based on this feature, they proposed a memetic algorithm which is very effective in solving UBQP problems with up to $2\,500$ variables. In the work of Tari et al. \cite{pivotingrules}, they used a recently proposed method, the monotonic local optima networks \cite{MLON}, to study the induced fitness landscapes for understanding the difference in performance among five pivoting rules.

One of the main challenge of combinatorial optimization comes from rugged landscapes. Some studies have proposed landscape smoothing based algorithms to improve algorithm performance. To the best of our knowledge, there are few landscape smoothing methods for UBQP, so we briefly introduce some related methods for other combinatorial optimization problems below. Gu and Huang \cite{GH} proposed a landscape smoothing method by edge cost manipulation. Liang et al. \cite{Liang} considered using local search to smooth the fitness landscape and thus preventing sticking in local optima for a long time. One year later, Coy et al. \cite{Coy} proposed the sequential smoothing algorithm (SSA) that alternates between convex and concave smoothing function to avoid being trapped in poor local optima. Hasegawa and Hiramatsu \cite{Hasegawa} suggested that metropolis algorithm can be used effectively as a local search algorithm in search-space smoothing strategies. Recently, Shi et al. \cite{HC} proposed a new landscape smoothing method called HC transformation for the TSP, and this approach was recently extended to UBQP \cite{HCUBQP}.

\section{Methods}\label{3}
\begin{figure}[htbp]
    \centering
    \includegraphics[scale=0.3]{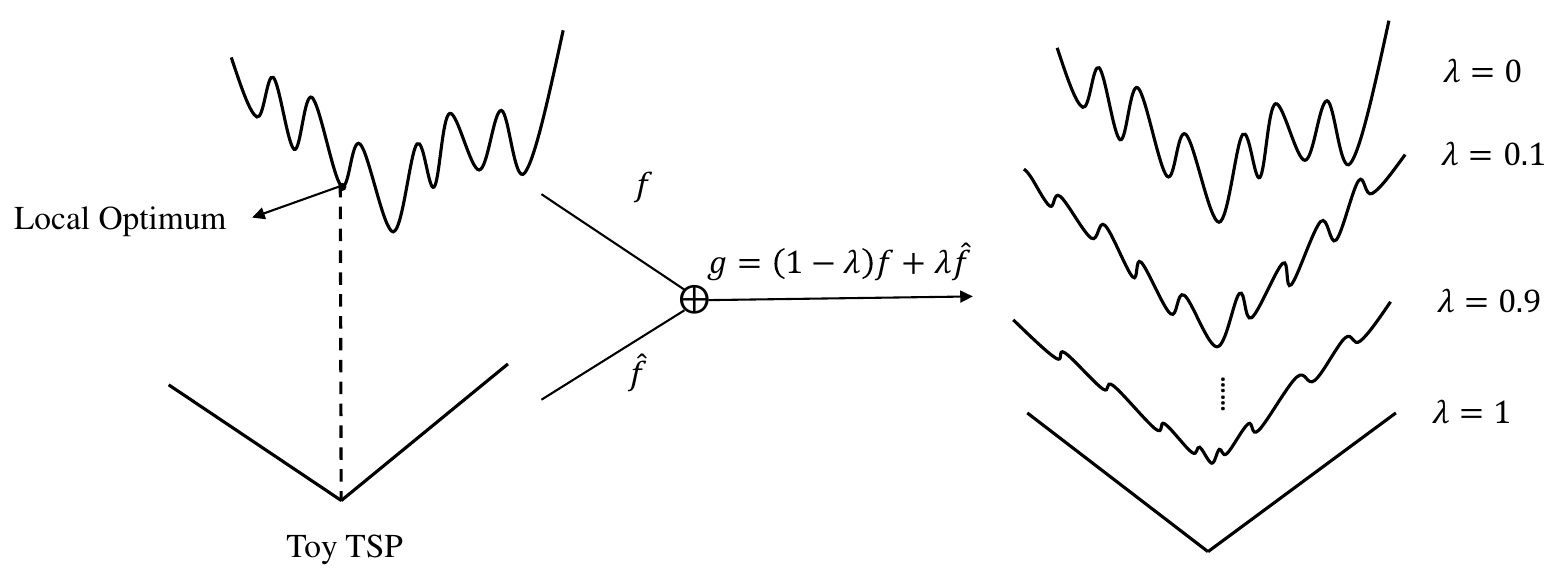}
    \caption{Effect of HC transformation for the TSP~\cite{HCUBQP}.}
    \label{Fig1}
\end{figure}

Initially proposed in \cite{HC}, HC transformation is a technique to smooth the rugged landscape of the traveling salesperson problem (TSP), with the detailed process illustrated in \figurename~\ref{Fig1}. The TSP can be smoothed by a convex combination of the original TSP and the unimodal TSP which is generated based on a known local optimum of the original TSP. This method was later generalized to UBQP~\cite{HCUBQP}. Given a known locally optimal solution for the considered UBQP instance, a toy UBQP of same size is created, represented by a binary matrix~$\boldsymbol{Q}$ constructed with $-1$s and $1$s values. Hence, we say that the toy UBQP used in \cite{HCUBQP} is constructed by the ``+/-1'' method. The unimodal nature of the toy UBQP landscape is theoretically established in \cite{HCUBQP}, ensuring that any local search process consistently converges towards its single global optimum, mirroring the known local optimum of the original UBQP. Subsequently, a convex combination of this constructed toy UBQP and the original UBQP is implemented, introducing a smoothing effect on the fitness landscape of the original UBQP. Governed by a coefficient denoted as $\lambda$ within the [0, 1] range, this process establishes a path from the original UBQP to the unimodal toy UBQP, allowing for an exploration of landscape characteristics and the impact of smoothing on search performance.

The LSILS algorithm represents an optimization approach that integrates the HC transformation technique with ILS. The HC transformation serves as a main component of LSILS, contributing to the landscape smoothing of the target UBQP instance. At the same time, ILS is used to explore and exploit the smoothed landscape, while concurrently updating the smoothed UBQP landscape iteratively to maintain the best solution found so far in the original UBQP. The entire procedure of LSILS is depicted in Algorithm \ref{LSILS}, where the algorithm executes an iterative process of a local search procedure and a perturbation procedure. This iterative process ensures the smoothing effect of the solution space and the enhancement of optimization outcomes.
\begin{algorithm}[t]
		\caption{Landscape Smoothing Iterated Local Search for UBQPs}
		\label{LSILS}
		$\boldsymbol{x}_{ini}\leftarrow$ Randomly generated solution\;
		$\boldsymbol{x}_{(0)}\leftarrow$ ILS$(\boldsymbol{x}_{ini}|f_o)$\;
		$\boldsymbol{x}^\star\leftarrow \boldsymbol{x}_{(0)}$\;
		$j\leftarrow 0$\;
		\While{stopping criterion is not met}
		{
			Construct the unimodal UBQP $\hat{f}$ based on $\boldsymbol{x}^\star$\;
			$g\leftarrow(1-\lambda)f_o+\lambda \hat{f}$\;
			$\boldsymbol{x}'_{(j)}\leftarrow$Perturbation($\boldsymbol{x}_{(j)}$)\;
			$\{\boldsymbol{x}_{(j+1)},\boldsymbol{x}^\star\}\leftarrow$LS($\boldsymbol{x}'_{(j)},\boldsymbol{x}^\star|g$)\;
			$j\leftarrow j+1$\;
			$\lambda$ = Update($\lambda$)\;
		}
		\Return $\boldsymbol{x}^\star$
\end{algorithm}

Lines 1-4 initiate the LSILS algorithm, initializing the current best solution denoted as $\boldsymbol{x}^\star$ with respect to $f_o$, where $f_o$ represents the original objective function. Subsequently, LSILS iteratively executes the ensuing steps until the termination criterion is satisfied (lines 5-11). At each iteration, the procedure starts by constructing the objective function $g$ for the smoothed UBQP based on the current best solution $\boldsymbol{x}^\star$ (line 7). This objective function $g$ is obtained by integrating the original UBQP objective function $f_o$ and the objective function $\hat{f}$ of the toy UBQP, employing a smoothing factor $\lambda$. After that, a local search operation is applied to the smoothed UBQP, initiating from a perturbed solution $\boldsymbol{x}'_{(j)}$ derived from the current solution $\boldsymbol{x}_{(j)}$ (line 8). Here, \mbox{LS($\boldsymbol{x}'_{(j)}$, $\boldsymbol{x}^\star|g$)} denotes the execution of a local search from $\boldsymbol{x}'_{(j)}$ on the smoothed UBQP $g$, concurrently updating $\boldsymbol{x}^\star$ by tracking the original objective function $f_o$. The resultant solution, denoted as $\boldsymbol{x}_{(j+1)}$, represents the local optimum for the smoothed UBQP~$g$ (line 9) and is utilized as the input for the next local search after perturbation (line 10). The final step involves the update of the smoothing coefficient $\lambda$ (line 11), with the specific methodology for this update being user-defined.

In~\cite{HCUBQP}, the construction of the toy UBQP (line 6 in Algorithm \ref{LSILS}) is based on a simple method called ``+/-1'', which is detailed as follows. Given a solution vector $\boldsymbol{x}^\star$ of dimension $n$, the generation of the specialized matrix $\boldsymbol{\hat{Q}}^1$ is executed through the following process,
\begin{equation}
	\label{Q1}
	\hat{Q}^1_{ij} = \left\{
			\begin{aligned}
				1,\quad
				\text{if}  \hspace{0.5em} x^\star_ix^\star_j=1\\
				-1,\quad \text{if}   \hspace{0.5em} x^\star_ix^\star_j\neq1
			\end{aligned}
		\right.,\ \ \  i,j\in\{1,2,\cdots,n\}
\end{equation}
where the $i$-th and $j$-th binary variables (0-1) of solution $\boldsymbol{x}^\star$ are denoted by $x^\star_i$ and $x^\star_j$, respectively. $\boldsymbol{\hat{Q}}^1 = [\hat{Q}^1_{ij}]$ represents an $n \times n$ matrix of the constructed toy UBQP. From the construction process of $\boldsymbol{\hat{Q}}^1$, we can see that only when the solution of the toy problem is $\boldsymbol{x}^\star$, all the positive elements in $\boldsymbol{\hat{Q}}^1$ are obtained, while for any solution except~$\boldsymbol{x}^\star$, a certain number of negative elements will be taken. This reflects the unimodality of the toy UBQP. In \cite{HCUBQP}, it is theoretically proven that the toy UBQP $\hat{f}^1$ defined by $\boldsymbol{\hat{Q}}^1$ exhibits unimodal behavior, indicating that a local search algorithm will consistently converge towards the unique global optimum of $\hat{f}^1$, i.e., $\boldsymbol{x}^\star$. The unimodal characteristic of $\hat{f}^1$ can be used to smooth the landscape of the original UBQP while preserving valuable information in the high quality solution $\boldsymbol{x}^\star$, thereby facilitating a more efficient search towards optimal solutions.

Given the unique features of the problem landscape, the efficiency of an algorithm is influenced not only by the smoothness of the landscape but also by its flatness. While smoothness denotes the global regularity within the solution space, flatness focuses on the horizontal characteristics of local regions. In this paper, in order to further examine the impact of the HC transformation on smoothing rugged landscapes by employing different toy problems, we offer two new constructions of the toy UBQP with different flatness, designated as $\boldsymbol{\hat{Q}}^2$ and $\boldsymbol{\hat{Q}}^3$.  

Similar to $\boldsymbol{\hat{Q}}^1$, the matrix $\boldsymbol{\hat{Q}}^2$ is constructed based on a method called ``+/-i'' as follows:
\begin{equation}
	\label{Q2}
	\hat{Q}^2_{ij} = \left\{
			\begin{aligned}
				i,\quad
				\text{if}  \hspace{0.5em} x^\star_ix^\star_j=1\\
				-i,\quad \text{if}   \hspace{0.5em} x^\star_ix^\star_j\neq1
			\end{aligned}
		\right.
\end{equation}
where $i\in\{1,2,\cdots,n\}$ and $j\in\{1,2,\cdots,i\}$, $i$ is the $i$-th row of matrix $\boldsymbol{\hat{Q}}^2$. Let $\hat{Q}^2_{ji} = \hat{Q}^2_{ij}$ to construct the symmetric $\boldsymbol{\hat{Q}}^2$. By contrast, the matrix $\boldsymbol{\hat{Q}}^3$ is constructed based on a so-called ``randomly'' method:
\begin{equation}
	\label{Q3}
	\hat{Q}^3_{ij} = \left\{
			\begin{aligned}
				u^1_{ij},\quad
				\text{if}  \hspace{0.5em} x^\star_ix^\star_j=1\\
				u^2_{ij},\quad \text{if}   \hspace{0.5em} x^\star_ix^\star_j\neq 1
			\end{aligned}
		\right.
\end{equation}
where $i\in\{1,2,\cdots,n\}$ and $j\in\{1,2,\cdots,i\}$, $u^1_{ij}$ is a randomly generated integer from a discrete uniform distribution in the range [1, 100], that is $P(u^1_{ij} = k_1) = \frac{1}{100}, k_1\in\{1,2,...,100\}$ and $u^2_{ij}$ is a randomly generated integer from a discrete uniform distribution in the range [-100, -1], that is $P(u^2_{ij} = k_2) = \frac{1}{100}, k_2 \in\{-100,-99,...,-1\}$. Then let $\hat{Q}^3_{ji} = \hat{Q}^3_{ij}$ to construct the symmetric $\boldsymbol{\hat{Q}}^3$. 

Since $\boldsymbol{\hat{Q}}^2$ (Eq.~(\ref{Q2})) and $\boldsymbol{\hat{Q}}^3$ (Eq.~(\ref{Q3})) have the same positive/negative element distribution as $\boldsymbol{\hat{Q}}^1$ (Eq.~(\ref{Q1})), which is determined by $\boldsymbol{x}^\star$, we can say that the landscapes defined by $\boldsymbol{\hat{Q}}^2$ and $\boldsymbol{\hat{Q}}^3$ are both unimodal with $\boldsymbol{x}^\star$ as the unique optimum.

As an example, for a 5-dimensional UBQP, given $\boldsymbol{x}^\star = (0,1,0,1,1)$, we construct three toy UBQPs based on the above constructions. Initially, it can be seen that $x^\star_1=0$, $x^\star_2=1$, $x^\star_3=0$, $x^\star_4=1$, $x^\star_5=1$, according to Eq.~(\ref{Q1}), Eq.~(\ref{Q2}) and Eq.~(\ref{Q3}), we have $\hat{Q}^1_{22} , \hat{Q}^1_{44} , \hat{Q}^1_{55} , \hat{Q}^1_{24} , \hat{Q}^1_{42} , \hat{Q}^1_{25} , \hat{Q}^1_{52} , \hat{Q}^1_{45} , \hat{Q}^1_{54} > 0 $, the remaining elements are all negative.
\[
\boldsymbol{\hat{Q}}^1 = \begin{pmatrix}
	-1 & -1 & -1 & -1 & -1 \\
	-1 & 1 & -1 & 1 & 1 \\
	-1 & -1 & -1 & -1 & -1 \\
	-1 & 1 & -1 & 1 & 1 \\
	-1 & 1 & -1 & 1 & 1 \\
\end{pmatrix}
\]

\[
\boldsymbol{\hat{Q}}^2 = \begin{pmatrix}
	-1 & -2 & -3 & -4 & -5 \\
	-2 & 2 & -3 & 4 & 5 \\
	-3 & -3 & -3 & -4 & -5 \\
	-4 & 4 & -4 & 4 & 5 \\
	-5 & 5 & -5 & 5 & 5 \\
\end{pmatrix}
\]

\[
\boldsymbol{\hat{Q}}^3 = \begin{pmatrix}
	-33 & -27 & -52 & -46 & -40 \\
	-27 & 62 & -72 & 95 & 11 \\
	-52 & -72 & -24 & -18 & -44 \\
	-46 & 95 & -18 & 1 & 17 \\
	-40 & 11 & -44 & 17 & 21 \\
\end{pmatrix}
\]

The aforementioned toy UBQP instances defined by $\boldsymbol{\hat{Q}}^m,~m\in\{1,2,3\}$ can be used to smooth the original UBQP. After smoothing, the objective function ($g$) is defined as:
\begin{equation} \nonumber
    \begin{aligned}
	   g(\boldsymbol{x}|\lambda,\alpha,m) &= (1-\lambda)f_o(\boldsymbol{x}) + \lambda \alpha  \hat{f}^m(\boldsymbol{x})\\
            &= (1-\lambda)\boldsymbol{x}^T\boldsymbol{Q}\boldsymbol{x} + \lambda \alpha \boldsymbol{x}^T\boldsymbol{\hat{Q}}^m\boldsymbol{x},~~m\in\{1,2,3\},
	\end{aligned}
\end{equation}
where $f_o$ is the objective function of the original UBQP, $\boldsymbol{Q}$ is the matrix of $f_o$,  $\hat{f}^m$ is the objective function of the toy UBQP defined by matrix $\boldsymbol{\hat{Q}}^m,~m\in\{1,2,3\}$. $\lambda\in[0,1]$ controls the strength of smoothing and $\alpha$ is a scaling factor to make sure that the original UBQP $f_o$ and the toy UBQP $\hat{f}^m$ are on the same scale. We can see that, when $\lambda=0$, the smoothed UBQP~$g$ degenerates to the original UBQP $f_o$, when $\lambda=1$, it is smoothed to the unimodal UBQP $\hat{f}^m$ and when \mbox{$0<\lambda<1$}, the original UBQP $f_o$ can be gradually smoothed to the unimodal UBQP $\hat{f}^m$. In fact, The smoothed objective function $g$ is a homotopic transformation from $f_o$ to $\hat{f}^m$.

\section{Experimental Analysis}\label{4}
To comprehensively investigate the impact of the HC transformation on smoothing rugged landscapes using various toy problems, a series of experiments are conducted. To better understand the solution distributions of different toy UBQPs, we first conduct landscape flatness experiments, and then we apply LSILS algorithms to the smoothed UBQPs, which have been modified by means of the HC transformation.

In order to carry out the experiments, we select 10 UBQP instances with a size of $n = 2\,500$ from the ORLIB \cite{orlib}. These instances possess a density of 0.1, where density denotes the proportion of non-zero elements in the matrix $\boldsymbol{Q}$. All 10 instances have known optimal solutions. The experiments are conducted using GNU C++ with the -O2 optimization option on the Tianhe-2 supercomputer, which is equipped with $17\,920$ computer nodes, each comprising two Intel Xeon E5-2692 12C (2.200 GHz) processors. 

The parameter setting is as follows: the maximum value of~$\lambda$ is set to 0.004. The time point to increase $\lambda$ is based on the CPU time, changing from 0 to 0.001 at 200 s, from 0.001 to 0.002 at 400 s, and so forth. This is because the smoothing effect of the HC transformation depends highly on the quality of the local optimum used to construct the toy UBQP. At the beginning of the search, the quality of the solution is poor, so we set a small $\lambda$. As search goes by, better solutions are found, a larger $\lambda$ is set to make full use of the high-quality solutions. Considering the zero elements of the original matrix~$\boldsymbol{Q}$ and the non-sparse nature of the toy UBQPs, to smooth the landscape of the original UBQP without causing great damage, we introduce $\alpha$ to constrain the elements of~$\boldsymbol{\hat{Q}}^m$ to be less than or equal to 5. Notably, the choice of 5 is proximate to the average absolute value of the original matrix~$\boldsymbol{Q}$. Specifically, for the three different toy UBQPs, $\alpha$ is set to be 5 for $\boldsymbol{\hat{Q}}^1$ (construct by ``+/-1''), 0.002 for~$\boldsymbol{\hat{Q}}^2$ (construct by ``+/-i'') and 0.05 for $\boldsymbol{\hat{Q}}^3$ (construct randomly). The stopping criterion is at $1\,000$ seconds of CPU runtime, with information logging conducted every 10 seconds. Each algorithm was executed for 20 independent runs on each benchmark instance, employing different random initial solutions. The local search implementation incorporates the \emph{1-bit-flip} local search, while perturbation is executed randomly, flipping $n/4$ bits at each call, where $n$ is the instance size.

\subsection{Landscape Analysis of Toy UBQPs}
Before the implementation of LSILS on UBQP instances, an initial experimental analysis is conducted to evaluate the landscape flatness of the three different toy problems. This comprehensive experiment aims to provide a thorough understanding of the landscape characteristics.

As a first step of the experiment, we use a 18-dimensional UBQP problem as an example. A 18-dimensional symmetric matrix is generated randomly as an original UBQP problem, whose elements follow a discrete uniform distribution in the range [-100,100]. Based on a solution $\boldsymbol{x}^\star=(1, 1, 0, 0, 1, 0, 0, 0, 0, 0, 1, 1, 1, 1, 1, 1, 1, 0)$, we build the corresponding $\boldsymbol{\hat{Q}}^1$ (construct by ``+/-1''), $\boldsymbol{\hat{Q}}^2$ (construct by ``+/-i'') and $\boldsymbol{\hat{Q}}^3$ (construct randomly) for the 18-dimensional randomly generated UBQP. Since the instance size is small, we traverse all solutions from the search space. 

\begin{figure}[htbp]
  \centering
  \includegraphics[width = 2.5in]{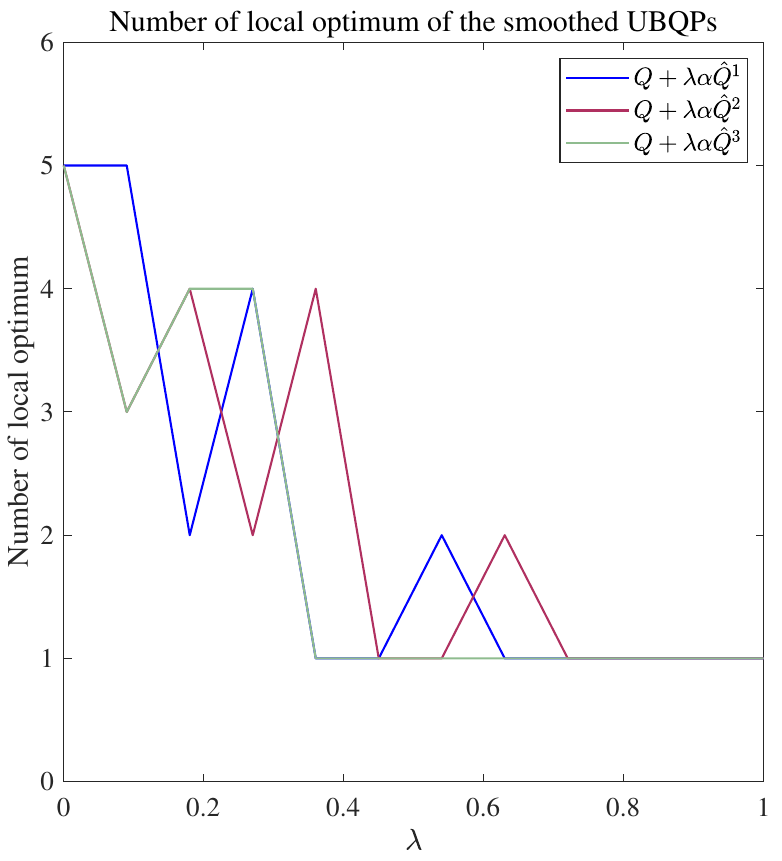}
  \caption{Number of local optima of the smoothed UBQPs with different $\lambda$ values.}\label{lambda}
\end{figure}

We first focus on the number of 1-bit-flip local optimum in the solution space of the four UBQPs which can reflect the ruggedness of the landscape. As shown in \tablename~\ref{localoptimanumber}, the original randomly generated $\boldsymbol{Q}$ has 5 1-bit-flip local optima, while $\boldsymbol{\hat{Q}}^1$ (construct by ``+/-1''), $\boldsymbol{\hat{Q}}^2$ (construct by ``+/-i'') and~$\boldsymbol{\hat{Q}}^3$ (construct randomly) all have only 1, that is $\boldsymbol{x}^\star$. This experiment confirms the unimodality of the three toy UBQPs. We also conduct HC transformation using the three toy UBQPs on the original randomly generated problem with different $\lambda$ values. To make sure the original UBQP and the three toy UBQPs are in the same scale, we let $\alpha = 50$ in the HC transformation with $\boldsymbol{\hat{Q}}^1$ (construct by  ``+/-1''), $\alpha = 2.8$ in the HC transformation with $\boldsymbol{\hat{Q}}^2$ (construct by  ``+/-i''), $\alpha = 1$ in the HC transformation with~$\boldsymbol{\hat{Q}}^3$ (construct randomly), where 50 is proximate to the average absolute value of the original randomly generated matrix $\boldsymbol{Q}$. \figurename~\ref{lambda} shows how the number of local optimum of the smoothed UBQP landscape changes against $\lambda$. In general, we can see that the number of local optima is approximately negatively related to $\lambda$ and when $\lambda > 0.8$, all the smoothed UBQPs tend to have only 1 local optimum. Based on the above observations, we can conclude that $\boldsymbol{\hat{Q}}^1$ (construct by ``+/-1''), $\boldsymbol{\hat{Q}}^2$ (construct by ``+/-i'') and~$\boldsymbol{\hat{Q}}^3$ (construct randomly) are unimodal and the HC transformation can indeed smooth the UBQP landscape, where the smoothing effect is controlled by $\lambda$. 

\begin{table}[htbp]
\centering
\caption{Local optima of four different 18-dimensional UBQPs.}
\label{localoptimanumber}
\begin{tabular}{|c|c|c|}
\hline
Instance                                                                 & \begin{tabular}[c]{@{}c@{}}Number of \\ local optimum\end{tabular} & 1-bit local optimum                                                                           \\ \hline
\begin{tabular}[c]{@{}c@{}}Original randomly \\ generated $\boldsymbol{Q}$\end{tabular} & 5                                                                 & \begin{tabular}[c]{@{}c@{}}$\boldsymbol{x}_1$ = (0, 1, 0, 0, 0, 0, 0, 1, 1, \\ 1, 0, 1, 1, 0, 0, 0, 0, 0)\\ $\boldsymbol{x}_2$ = (1, 0, 1, 1, 1, 0, 0, 0, 1, \\ 0, 1, 1, 1, 0, 0, 0, 1, 0)\\ $\boldsymbol{x}_3$ = (1, 1, 1, 0, 1, 0, 0, 1, 1, \\ 0, 1, 1, 1, 0, 0, 0, 1, 0)\\ $\boldsymbol{x}_4$ = (0, 0, 0, 1, 1, 0, 0, 1, 1, \\ 1, 1, 1, 1, 0, 0, 0, 1, 0)\\ $\boldsymbol{x}_5$ = (1, 1, 1, 0, 1, 0, 0, 1, 1, \\ 0, 0, 1, 1, 0, 1, 0, 1, 0)\end{tabular} \\ \hline
\begin{tabular}[c]{@{}c@{}}$\boldsymbol{\hat{Q}}^1$ (construct \\ by  ``+/-1'')\end{tabular}      & 1                                                                 & \begin{tabular}[c]{@{}c@{}}$\boldsymbol{x}_1$ = (1, 1, 0, 0, 1, 0, 0, 0, 0, \\ 0, 1, 1, 1, 1, 1, 1, 1, 0)\end{tabular}                                                                                                                                                            \\ \hline
\begin{tabular}[c]{@{}c@{}}$\boldsymbol{\hat{Q}}^2$ (construct \\ by ``+/-i'')\end{tabular}      & 1                                                                 & \begin{tabular}[c]{@{}c@{}}$\boldsymbol{x}_1$ = (1, 1, 0, 0, 1, 0, 0, 0, 0, \\ 0, 1, 1, 1, 1, 1, 1, 1, 0)\end{tabular}                                                                                                     \\ \hline
\begin{tabular}[c]{@{}c@{}}$\boldsymbol{\hat{Q}}^3$ (construct \\ randomly)\end{tabular}       & 1                                                                 & \begin{tabular}[c]{@{}c@{}}$\boldsymbol{x}_1$ = (1, 1, 0, 0, 1, 0, 0, 0, 0, \\ 0, 1, 1, 1, 1, 1, 1, 1, 0)\end{tabular}                                                                                                     \\ \hline
\end{tabular}
\end{table}

\begin{figure*}[htbp]
    \centering
	\subfloat[Original randomly generated $\boldsymbol{Q}$
]{\includegraphics[width=3.5in]{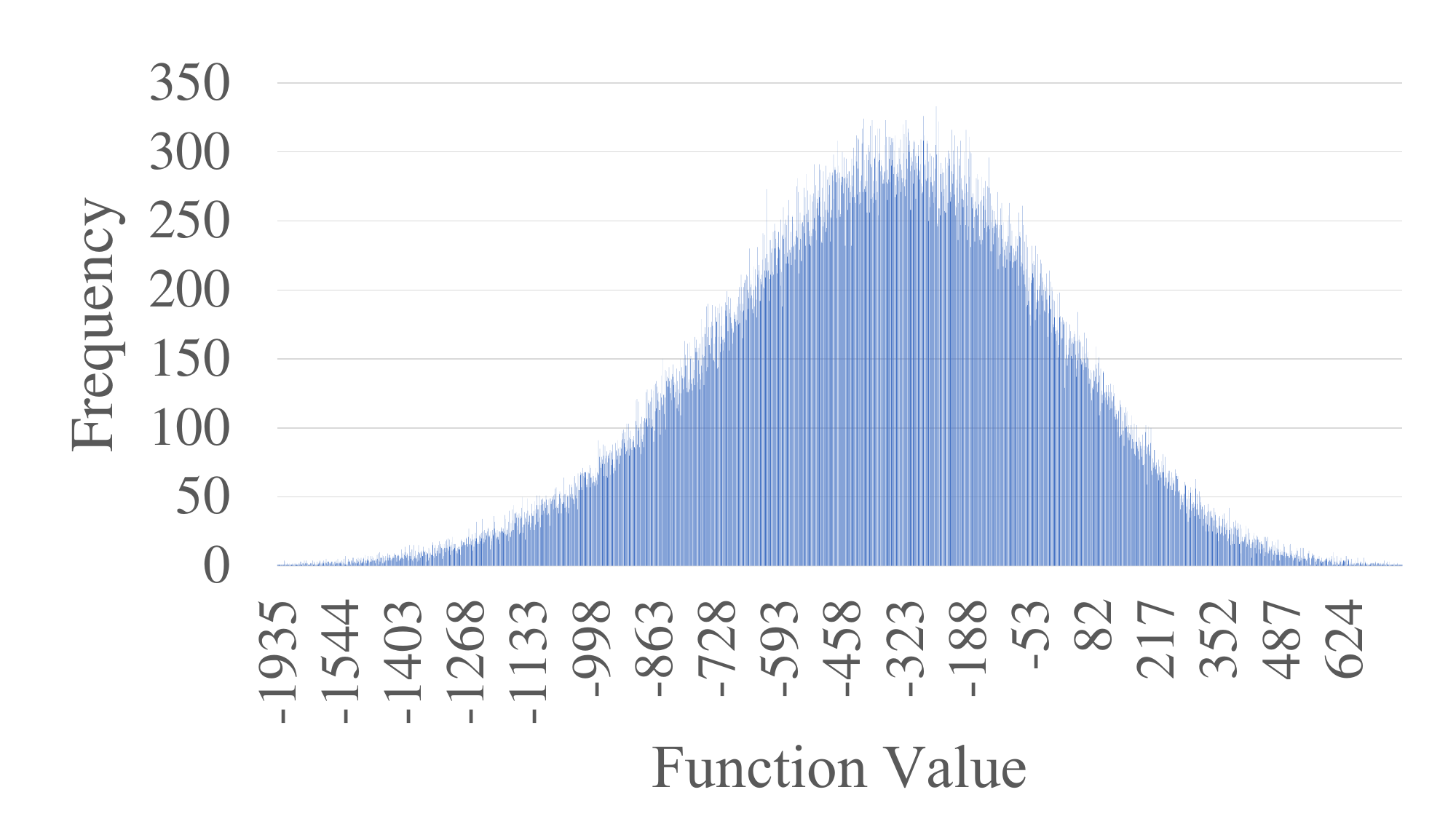}
		\label{original}}
	\subfloat[$\boldsymbol{\hat{Q}}^1$ (construct by ``+/-1'')]{\includegraphics[width=3.5in]{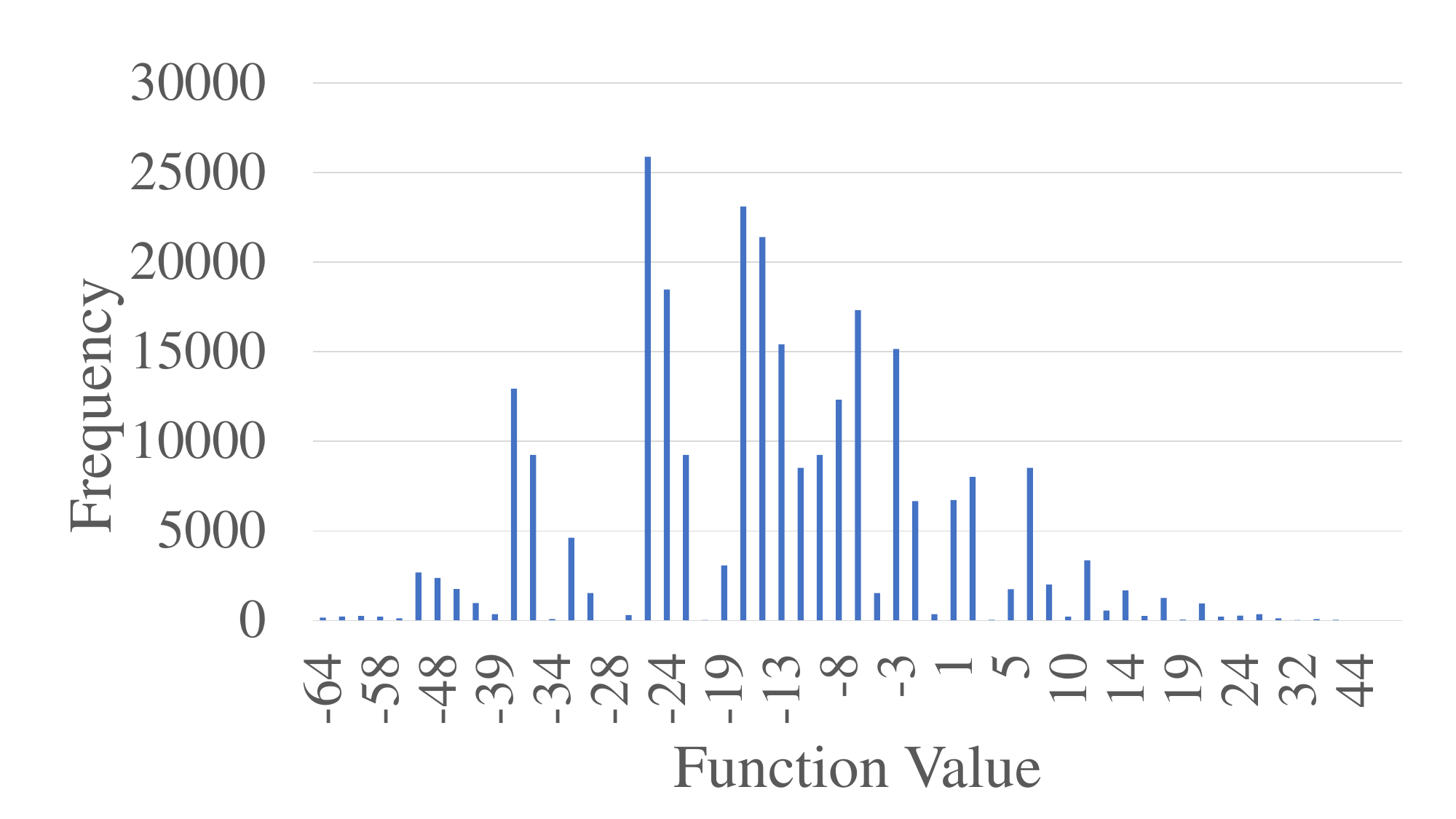}
		\label{1Q}}

	\subfloat[$\boldsymbol{\hat{Q}}^2$ (construct by ``+/-i'')]{\includegraphics[width=3.5in]{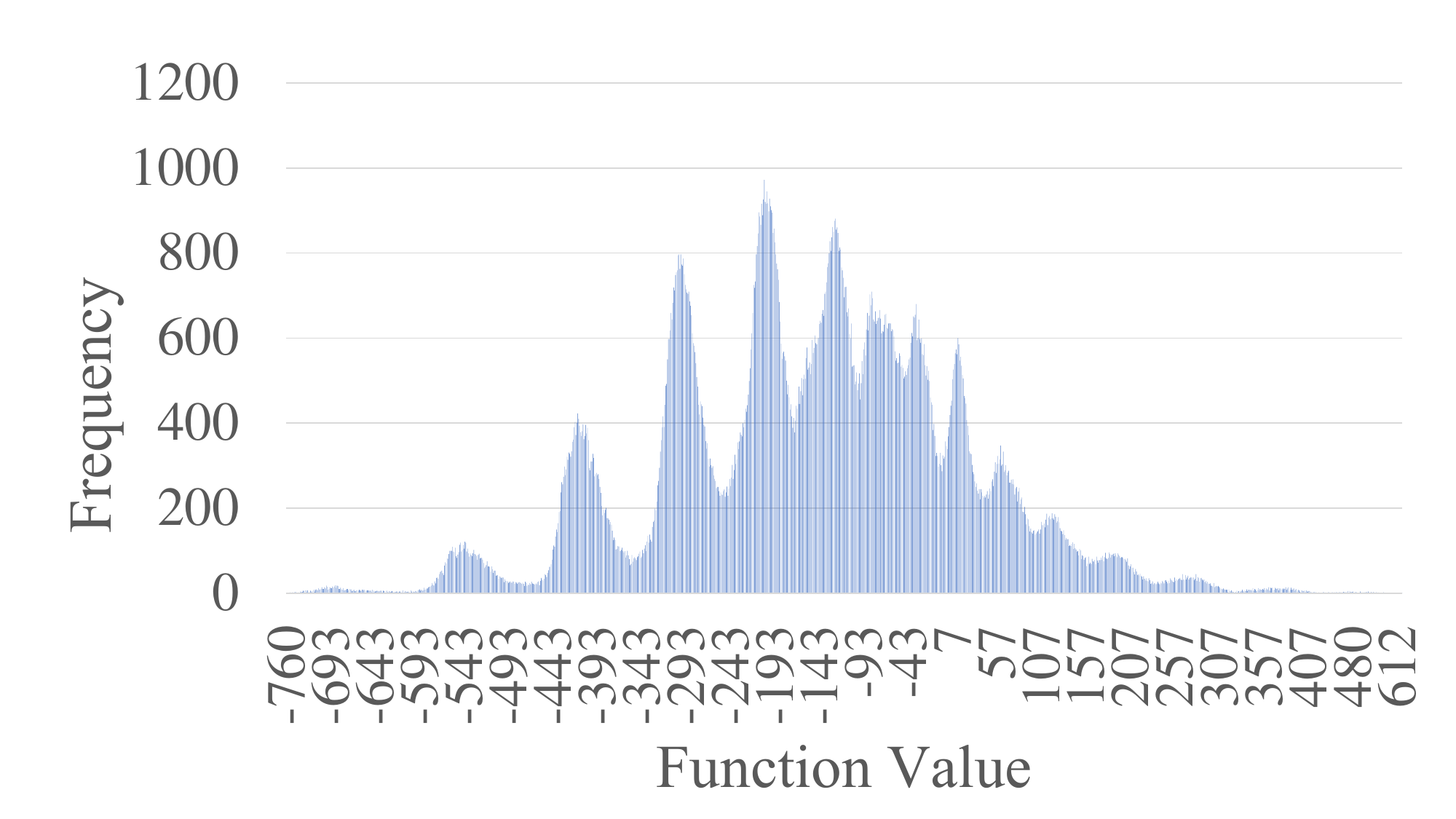}
		\label{iQ}}
	\subfloat[$\boldsymbol{\hat{Q}}^3$ (construct randomly)]{\includegraphics[width=3.5in]{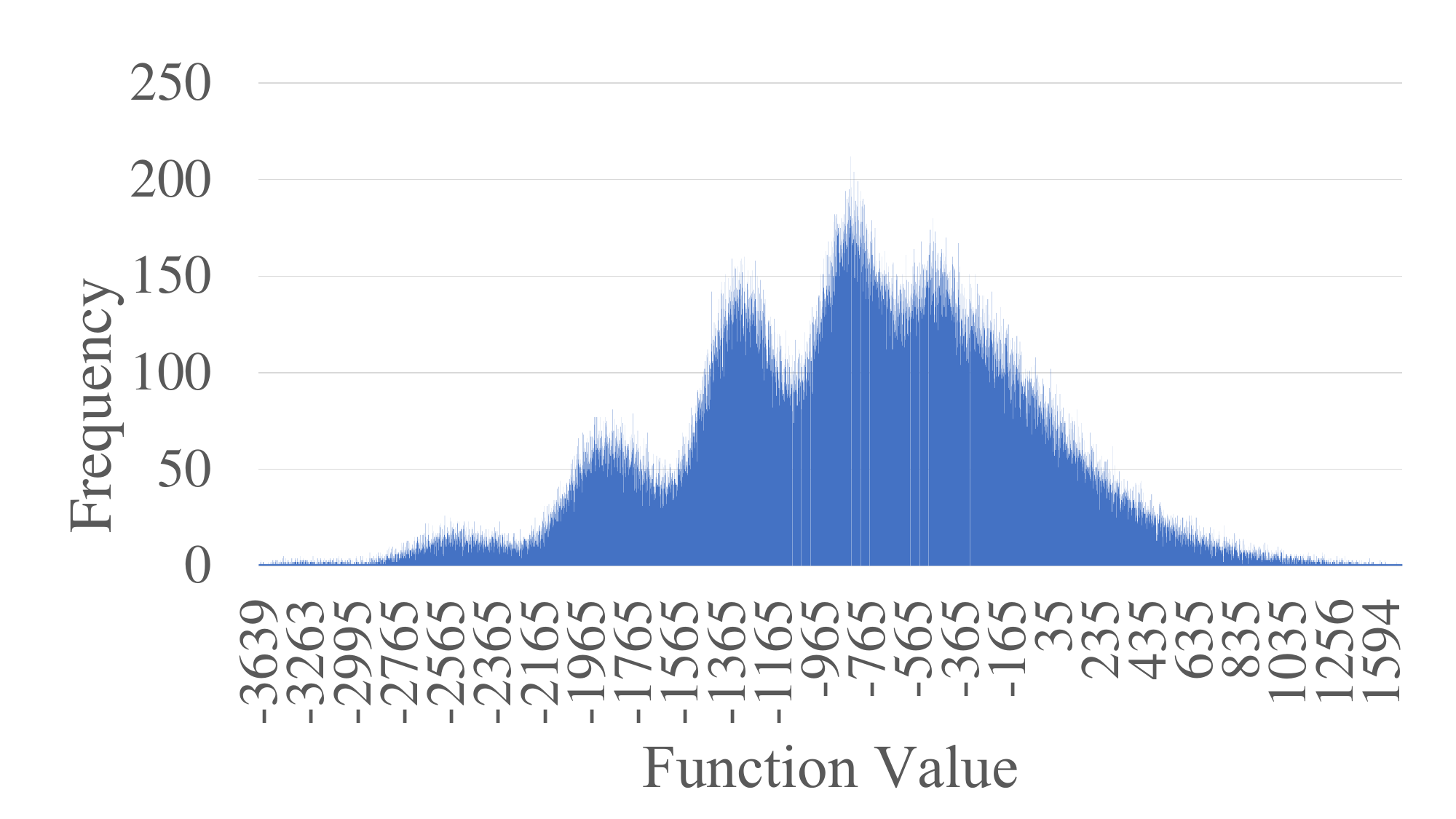}
		\label{randomQ}}
	\caption{Distribution of function values in the solution space of different 18-dimensional UBQPs, frequency versus the function value. (a) Original randomly generated $\boldsymbol{Q}$. (b) $\boldsymbol{\hat{Q}}^1$ (construct by ``+/-1''). (c) $\boldsymbol{\hat{Q}}^2$ (construct by ``+/-i''). (d) $\boldsymbol{\hat{Q}}^3$ (construct randomly).}
	\label{Distribution}
\end{figure*}

We then focus on the flatness of the original randomly generated $\boldsymbol{Q}$, $\boldsymbol{\hat{Q}}^1$ (construct by ``+/-1''), $\boldsymbol{\hat{Q}}^2$ (construct by ``+/-i'') and $\boldsymbol{\hat{Q}}^3$ (construct randomly). \figurename~\ref{original} gives a histogram of the original $\boldsymbol{Q}$ with function values as x axis and the corresponding frequency of each function value as y axis. From \figurename~\ref{original}, we can see that the function values tend to cluster around the central region and relatively few values are distributed at both ends, with a range of function values between -1935 and 922. To examine the distributions of function values of the three different toy UBQPs, an 18-dimensional vector is randomly generated from the solution space of the original UBQP, serving as the basis for constructing the toy UBQPs. By following the steps specified for constructing toy UBQPs in Eqs. (\ref{Q1}), (\ref{Q2}) and (\ref{Q3}), three different toy UBQP instances are obtained. The distributions of function values in the solution space for these three toy UBQPs are illustrated in \figurename~\ref{1Q}, \figurename~\ref{iQ}, and \figurename~\ref{randomQ}, respectively. Results reveal that the toy UBQP with $\boldsymbol{\hat{Q}}^1$ (construct by ``+/-1'') exhibits a relatively flat landscape, with more solutions sharing the same function values. In contrast, $\boldsymbol{\hat{Q}}^3$ (construct randomly) possesses a relatively non-flat landscape, with fewer solutions with the same function values. The landscape flatness of the toy UBQP with $\boldsymbol{\hat{Q}}^2$ (construct by ``+/-i'') falls between $\boldsymbol{\hat{Q}}^1$ (construct by ``+/-1'') and $\boldsymbol{\hat{Q}}^3$ (construct randomly), being flatter than $\boldsymbol{\hat{Q}}^3$ (construct randomly).

\begin{figure}[htbp]
    \centering
    \includegraphics[scale=0.55]{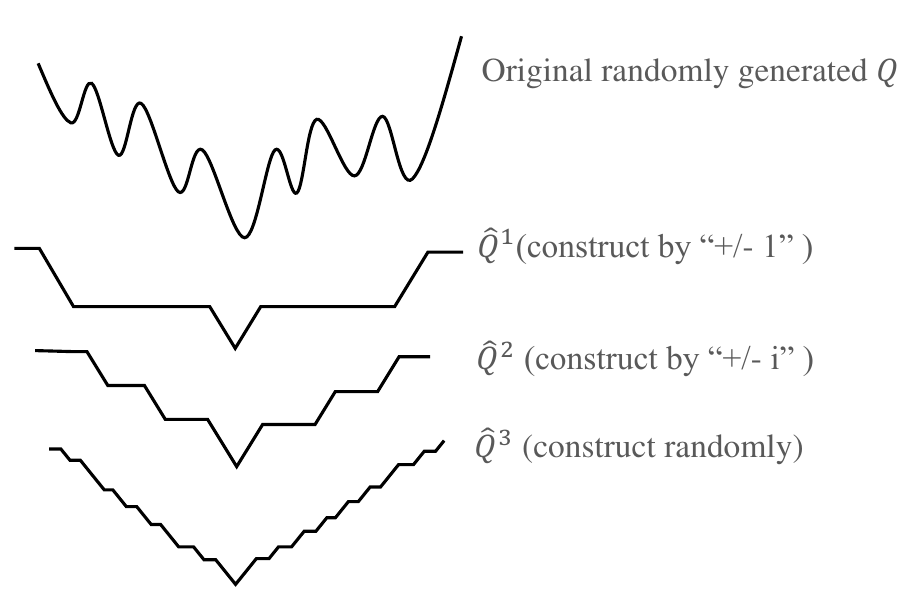}
    \caption{A diagram of the landscapes of different 18-dimensional UBQPs.}
    \label{Landscape}
\end{figure}

To quantify landscape flatness more precisely, the possibility of obtaining two solutions with different function values is calculated. This involves first summing the function value possibilities and then squaring the result. The resulting possibilities for the original UBQP, toy UBQP with $\boldsymbol{\hat{Q}}^1$ (construct by ``+/-1''), toy UBQP with $\boldsymbol{\hat{Q}}^2$ (construct by ``+/-i'), and toy UBQP with $\boldsymbol{\hat{Q}}^3$ (construct randomly) are 0.000795, 0.054152, 0.001924, and 0.000413, respectively. This metric reflects that though the landscape of the toy UBQP with $\boldsymbol{\hat{Q}}^3$ (construct randomly) is not as flat as the original UBQP, it remains unimodal. Conversely, the landscape of the toy UBQP with $\boldsymbol{\hat{Q}}^1$ (construct by ``+/-1'') appears excessively flat. The varying flatness of the toy UBQPs will imply diverse landscape smoothing effects on the original UBQP while owning the unimodal structure. The landscapes of the four UBQP problems are sketched in \figurename~\ref{Landscape}.

\subsection{Performance of LSILS with Different Toy UBQPs}
According to the previous analysis, three different toy UBQPs we designed have different landscapes. To study the landscape smoothing effects of the HC transformation with different toy UBQPs, we then conduct experiments on 10 large UBQP instances using LSILS. To measure the quality of a solution, we use its relative deviation to the optimum, which is defined as:
\begin{equation}
	\text{excess} = \frac{f(\boldsymbol{x}_{\textrm{LO}})-f(\boldsymbol{x}_{\text{opt}})}{f(\boldsymbol{x}_{\text{opt}})},
\end{equation}
where $f(\boldsymbol{x}_{\text{opt}})$ denotes the known global optimum and $f(\boldsymbol{x}_{\textrm{LO}})$ is the solution returned by the algorithm. In order to compare the effect of LSILS algorithms with different toy UBQPs, we use the result of ILS and LSILS using the toy UBQP with $\boldsymbol{\hat{Q}}^1$ (construct by ``+/-1'') from \cite{HCUBQP}. ILS iteratively executes a local search procedure and a perturbation procedure on the original UBQP landscape until the stopping criterion is met. The experimental settings are the same as before.

\begin{figure*}
    \centering
    \includegraphics[scale=0.4]{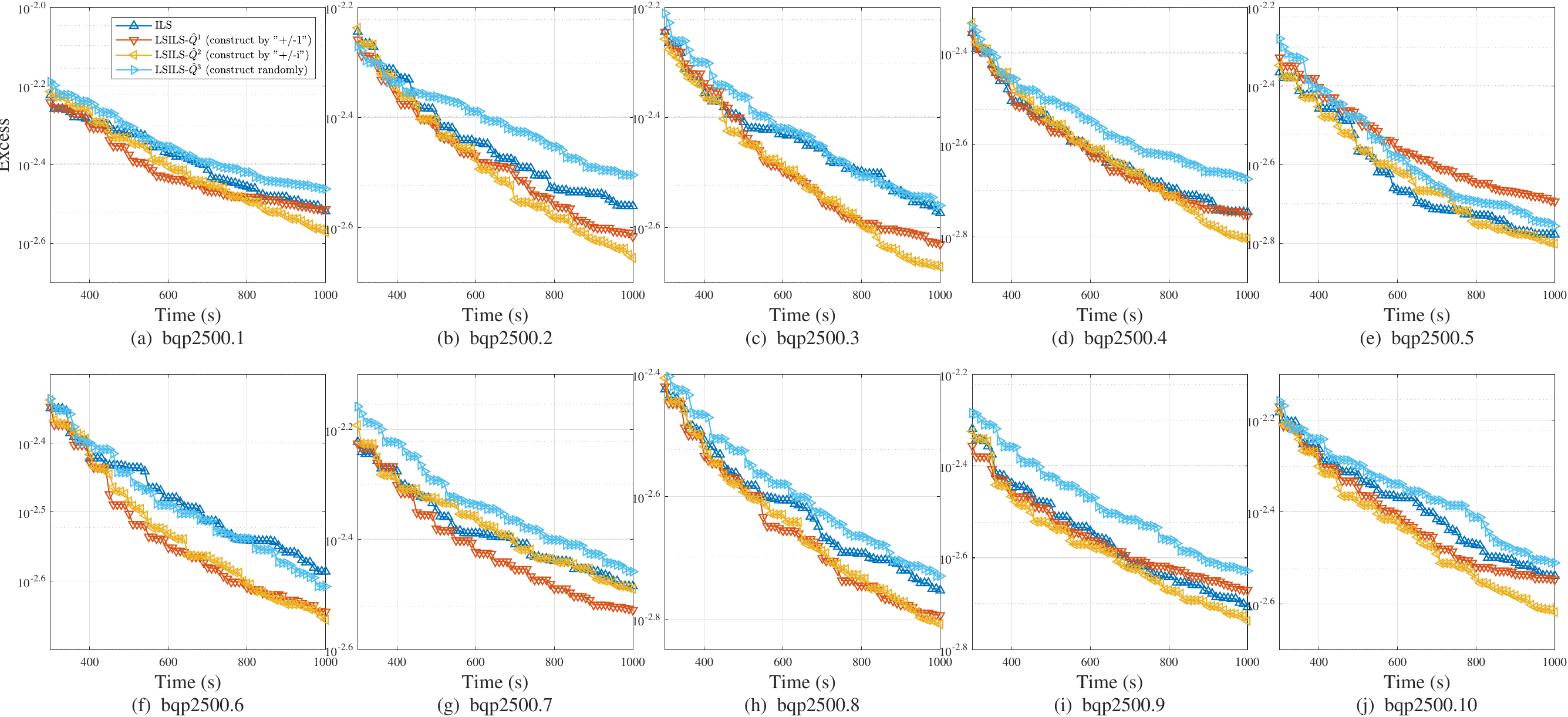}
    \caption{Comparison results of ILS, LSILS with different toy UBQPs. (a) bqp2500.1. (b) bqp2500.2. (c) bqp2500.3. (d) bqp2500.4. (e) bqp2500.5. (f) bqp2500.6. (g) bqp2500.7. (h) bqp2500.8. (i) bqp2500.9. (j) bqp2500.10.}
    \label{UBQP}
\end{figure*}

\figurename~\ref{UBQP} reports the average relative deviation of the best solution found by LSILS and ILS against CPU time. A lower value is better. From \figurename~\ref{UBQP}, we can see that in all instances, LSILS using the toy UBQP with $\boldsymbol{\hat{Q}}^2$ (construct by ``+/-i'') performs better than ILS. LSILS using the toy UBQP with $\boldsymbol{\hat{Q}}^1$ (construct by ``+/-1'') performs better than ILS on 8 instances except bqp2500.5 and bqp2500.9. In 9 out of 10 instances, LSILS using the toy UBQP with $\boldsymbol{\hat{Q}}^3$ (construct randomly) performs worse than ILS. These results compared with ILS confirm the effectiveness of LSILS using the toy UBQP with $\boldsymbol{\hat{Q}}^1$ (construct by ``+/-1'') and $\boldsymbol{\hat{Q}}^2$ (construct by ``+/-i'') on smoothing the landscape of the original UBQP. In contrast, the toy UBQP with $\boldsymbol{\hat{Q}}^3$ (construct randomly) has no landscape smoothing effect and even makes the original UBQP more difficult to solve. From previous landscape flatness analysis, we can find that the toy UBQP with $\boldsymbol{\hat{Q}}^3$ (construct randomly) possesses a non-flat landscape, with a lower likelihood of generating solutions with identical function values than the original UBQP.

In addition, \figurename~\ref{UBQP} shows that among the 10 instances, LSILS using the toy UBQP with $\boldsymbol{\hat{Q}}^2$ (construct by ``+/-i'') ultimately achieves better performance than with $\boldsymbol{\hat{Q}}^1$ (construct by ``+/-1''), except for the bqp2500.7 instance. This indicates that the toy UBQP with $\boldsymbol{\hat{Q}}^2$ (construct by ``+/-i'') has better landscape smoothing effect than with $\boldsymbol{\hat{Q}}^1$ (construct by ``+/-1''). Based on the landscape flatness analysis mentioned previously, a possible reason is that the toy UBQP with $\boldsymbol{\hat{Q}}^1$ (construct by ``+/-1'') has flatter landscape and these flat areas actually have almost no smoothing effect on the original UBQP by HC transformation, thus remaining the rugged landscape of the original UBQP landscape. These results reveal that the smoothing effect of the HC transformation is not only affected by the smoothness of the landscape but also by the flatness of the landscape, therefore, striking a balance between smoothness and flatness of the toy UBQP is crucial for HC transformation to guide the landscape smoothing algorithms.

In summary, the experimental results highlight that the toy UBQP with $\boldsymbol{\hat{Q}}^1$ (construct by ``+/-1'') displays the flattest landscape among the three toy UBQPs and the toy UBQP with $\boldsymbol{\hat{Q}}^3$ (construct randomly) demonstrates the most non-flat landscape while remaining unimodal. This landscape feature indeed makes the LSILS algorithm, which is based on the HC transformation, performs differently for different toy UBQPs. We can not simply assume that the flatter the landscape, the better the smoothing effect. To have a better landscape smoothing effect, choosing a suitable toy UBQP is a key issue.

	
	

\section{Conclusion}\label{5}
In this study, we extended the HC transformation based landscape smoothing method for UBQP by studying its impact when employing different toy problems: the toy UBQP with $\boldsymbol{\hat{Q}}^1$ (construct by ``+/-1''), the toy UBQP with $\boldsymbol{\hat{Q}}^2$ (construct by ``+/-i''), and the toy UBQP with $\boldsymbol{\hat{Q}}^3$ (construct randomly). Firstly, landscape analysis experiments were conducted to assess the flatness of the three toy UBQPs. The results demonstrate that the toy UBQP with $\boldsymbol{\hat{Q}}^1$ (construct by ``+/-1'') has the flattest landscape, while the toy UBQP with $\boldsymbol{\hat{Q}}^3$ (construct randomly)  exhibits the most non-flat landscape among the three. Secondly, the LSILS algorithm was applied to 10 UBQP benchmark instances utilizing each of the aforementioned toy UBQPs. Our findings show a different performance of LSILS with the three toy UBQPs, where LSILS using the toy UBQP with $\boldsymbol{\hat{Q}}^1$ (construct by ``+/-1'') and with $\boldsymbol{\hat{Q}}^2$ (construct by ``+/-i'') both having beneficial smoothing effects on the original UBQP. Among the two, the latter consistently outperforms the former. However, LSILS using the toy UBQP with $\boldsymbol{\hat{Q}}^3$ (construct randomly) even has a poorer performance than ILS. This reinforces the notion that different toy UBQPs exercise varying landscape smoothing effects on the original UBQP. The smoothness and flatness are both important landscape characteristics of the toy problem, emphasizing the important role of selecting an appropriate toy UBQP in enhancing the performance of LSILS. In conclusion, our study contributes valuable insights into the area of landscape smoothing through the utilization of the HC transformation with different toy problems. The observed variations in landscape flatness among the toy UBQPs directly influence the effectiveness of LSILS. In the future, we plan to investigate deep learning (DL) models in order to smooth difficult optimization problems whose solution spaces are rugged.

\section*{Acknowledgment}
The work described in this paper was supported in part by Sichuan Science and Technology Program (No. 2022NSFSC1907), in part by National Natural Science Foundation of China (grant no. 62076197, 61903294), in part by Key Research and Development Project of Shaanxi Province (grant no. 2022GXLH-01-15), in part by Research Grants Council of the Hong Kong Special Administrative Region, China (GRF Project No. CityU11215622), in part by Guangdong Basic and Applied Basic Research Foundation (No.2021A1515110077) and in part by Fundamental Research Founds for Central Universities (No.D5000210692).


\begin{thebibliography}{00}
\bibitem{UBQP} G. Kochenberger et al., “The unconstrained binary quadratic programming problem: a survey,” \textit{J Comb Optim}, vol. 28, no. 1, pp. 58–81, Jul. 2014.

\bibitem{financial} R. D. McBride and J. S. Yormark, “An Implicit Enumeration Algorithm for Quadratic Integer Programming,” \textit{Management Science}, vol. 26, no. 3, pp. 282–296, Mar. 1980.

\bibitem{molecular} A. T. Phillips and J. B. Rosen, “A quadratic assignment formulation of the molecular conformation problem,” \textit{J Glob Optim}, vol. 4, no. 2, pp. 229–241, Mar. 1994.

\bibitem{traffic} G. Gallo, P. L. Hammer, and B. Simeone, “Quadratic knapsack problems,” \textit{Combinatorial Optimization}, pp. 132–149. 1980.

\bibitem{UBQPformulation} A. Liefooghe, S. Verel, and J.-K. Hao, “A hybrid metaheuristic for multiobjective unconstrained binary quadratic programming,” \textit{Applied Soft Computing}, vol. 16, pp. 10–19, Mar. 2014.

\bibitem{maxcut} I. M. Bomze, M. Budinich, P. M. Pardalos, and M. Pelillo, “The Maximum Clique Problem,” in \textit{Handbook of Combinatorial Optimization: Supplement Volume A}, D.-Z. Du and P. M. Pardalos, Eds., Boston, MA: Springer US, 1999, pp. 1–74.

\bibitem{setpartition} M. Lewis, G. Kochenberger, and B. Alidaee, “A new modeling and solution approach for the set-partitioning problem,” \textit{Computers \& Operations Research}, vol. 35, no. 3, pp. 807–813, Mar. 2008.

\bibitem{otherproblem} G. A. Kochenberger, F. Glover, B. Alidaee, and C. Rego, “A unified modeling and solution framework for combinatorial optimization problems,” \textit{OR Spectrum}, vol. 26, no. 2, pp. 237–250, Mar. 2004.

\bibitem{np_hard} M. R. Garey, and D. S. Johnson, “Computers and Intractability: A Guide to the Theory of NP-Completeness,”. W. H. Freeman \& Co Ltd, 1979.

\bibitem{HCUBQP} W. Wang, J. Shi, J. Sun, A. Liefooghe, and Q. Zhang, “A New Parallel Cooperative Landscape Smoothing Algorithm and Its Applications on TSP and UBQP.” \textit{arXiv e-prints} arXiv:2401.03237, 2024.

\bibitem{ILS} H. R. Lourenço, O. C. Martin, and T. Stützle, “Iterated Local Search: Framework and Applications,” in \textit{Handbook of Metaheuristics}, M. Gendreau and J.-Y. Potvin, Eds., in International Series in Operations Research \& Management Science. , Cham: Springer International Publishing, 2019, pp. 129–168.

\bibitem{DDTS} F. Glover, Z. L$\ddot{u}$, and J.-K. Hao, “Diversification-driven tabu search for unconstrained binary quadratic problems,” \textit{4OR-Q J Oper Res}, vol. 8, no. 3, pp. 239–253, Oct. 2010.

\bibitem{TS} Z. L$\ddot{u}$, F. Glover, and J.-K. Hao, “A hybrid metaheuristic approach to solving the UBQP problem,” \textit{European Journal of Operational Research}, vol. 207, no. 3, pp. 1254–1262, Dec. 2010.

\bibitem{SA1} T. M. Alkhamis, M. Hasan, and M. A. Ahmed, “Simulated annealing for the unconstrained quadratic pseudo-Boolean function,” \textit{European Journal of Operational Research}, vol. 108, no. 3, pp. 641–652, Aug. 1998.

\bibitem{SA2} K. Katayama and H. Narihisa, “Performance of simulated annealing-based heuristic for the unconstrained binary quadratic programming problem,” \textit{European Journal of Operational Research}, vol. 134, no. 1, pp. 103–119, Oct. 2001.

\bibitem{ILS1} K. Katayama and H. Narihisa, “On fundamental design of parthenogenetic algorithm for the binary quadratic programming problem,” in \textit{Proceedings of the 2001 Congress on Evolutionary Computation}, 2001, pp. 356–363.

\bibitem{EA} A. Lodi, K. Allemand, and T. M. Liebling, “An evolutionary heuristic for quadratic 0–1 programming,” \textit{European Journal of Operational Research}, vol. 119, no. 3, pp. 662–670, Dec. 1999.

\bibitem{fitness_landscape} K. M. Malan, “A survey of advances in landscape analysis for optimisation,” \textit{Algorithms}, vol. 14, no. 2, pp. 40, 2021.

\bibitem{SampledWalk} S. Tari, M. Basseur, and A. Goëffon, “Sampled Walk and Binary Fitness Landscapes Exploration,” in \textit{Artificial Evolution}, E. Lutton, P. Legrand, P. Parrend, N. Monmarché, and M. Schoenauer, Eds., in Lecture Notes in Computer Science. Cham: Springer International Publishing, 2018, pp. 47–57.

\bibitem{UBQPworst} S. Tari, M. Basseur, and A. Goëffon, “Worst Improvement Based Iterated Local Search,” in \textit{Evolutionary Computation in Combinatorial Optimization}, A. Liefooghe and M. López-Ibáñez, Eds., in Lecture Notes in Computer Science. Cham: Springer International Publishing, 2018, pp. 50–66.

\bibitem{MA}P. Merz and K. Katayama, “Memetic algorithms for the unconstrained binary quadratic programming problem,” \textit{Biosystems}, vol. 78, no. 1, pp. 99–118, Dec. 2004.

\bibitem{pivotingrules} S. Tari and G. Ochoa, “Local search pivoting rules and the landscape global structure,” in \textit{Proceedings of the Genetic and Evolutionary Computation Conference, in GECCO ’21}, Jun. 2021, pp. 278–286.

\bibitem{MLON} G. Ochoa, N. Veerapen, F. Daolio, and M. Tomassini, “Understanding Phase Transitions with Local Optima Networks: Number Partitioning as a Case Study,” in \textit{Evolutionary Computation in Combinatorial Optimization}, B. Hu and M. López-Ibáñez, Eds., in Lecture Notes in Computer Science. Cham: Springer International Publishing, 2017, pp. 233–248.

\bibitem{GH} J. Gu and X. Huang, “Efficient local search with search space smoothing: a case study of the traveling salesman problem (TSP),” \textit{IEEE Trans. Syst., Man, Cybern.}, vol. 24, no. 5, pp. 728–735, May 1994.

\bibitem{Liang} K. H. Liang, X. Yao, and C. Newton, “Combining landscape approximation and local search in global optimization,” in \textit{Proceedings of the 1999 Congress on Evolutionary Computation-CEC99}, Washington, DC, USA: IEEE, 1999, pp. 1514–1520.

\bibitem{Coy} S. P. Coy, B. L. Golden, and E. A. Wasil, “A computational study of smoothing heuristics for the traveling salesman problem,” \textit{European Journal of Operational Research}, vol. 124, no. 1, pp. 15–27, Jul. 2000.

\bibitem{Hasegawa} M. Hasegawa and K. Hiramatsu, “Mutually beneficial relationship in optimization between search-space smoothing and stochastic search,” \textit{Physica A: Statistical Mechanics and its Applications}, vol. 392, no. 19, pp. 4491–4501, Oct. 2013.

\bibitem{HC}  J. Shi, J. Sun, Q. Zhang, and K. Ye, “Homotopic Convex Transformation: A New Landscape Smoothing Method for the Traveling Salesman Problem,” \textit{IEEE Trans. Cybern.}, vol. 52, no. 1, pp. 495–507, Jan. 2022.

\bibitem{orlib} J. E. Beasley, “Obtaining test problems via Internet,” \textit{J Glob Optim}, vol. 8, no. 4, pp. 429–433, Jun. 1996.
\end{thebibliography}
\end{document}